\documentclass{amsart}

\newtheorem{theorem}{Theorem}
\newtheorem{lemma}[theorem]{Lemma}

\theoremstyle{definition}
\newtheorem{definition}[theorem]{Definition}
\newtheorem{question}[theorem]{Question}
\theoremstyle{remark}
\newtheorem{remark}[theorem]{Remark}

\usepackage{empheq}
\usepackage{tikz}
\usepackage{amssymb}
\usepackage{stmaryrd}
\usepackage{enumerate}
\usepackage{multicol}
\usepackage[shortlabels]{enumitem}
\usepackage{calc}
\usepackage{float}

\usetikzlibrary{automata,positioning}
\usetikzlibrary{arrows}

\usepackage[utf8]{inputenc}

\usepackage[left=2.0cm,%
                right=2.0cm,%
                top=2.5cm,%
                bottom=2.5cm,%
                headheight=12pt,%
                a4paper]{geometry}%

\begin{document}

\title{A note on public announcements in standard G\"odel modal logic}

\author{Nicholas Pischke}
\address{Hoch-Weiseler Str. 46, Butzbach, 35510, Hesse, Germany}
\email{pischkenicholas@gmail.com}
\date{\today}

\keywords{Fuzzy Logic, Modal Logic, Public Announcements, G\"odel Logic}

\begin{abstract}
We study public announcement operators in the context of standard G\"odel modal logic as introduced by Caicedo and Rodriguez. Over that base logic, admitting a natural semantics over $[0,1]$-valued generalizations of modal Kripke models, we exhibit three possible semantic interpretations of the public announcement operator, all equivalent in a classical setting, and show that these logics are all different in terms of expressive strength. We provide partial completeness results for these logics via Hilbert-style calculi.
\end{abstract}

\maketitle

\section{Introduction and Preliminaries}
Public announcement operators for modal logics first occurred in \cite{Pla1989} and mark, in some sense, the beginning of the whole subject of dynamic epistemic logic (see \cite{DHK2008} for a textbook reference). 

These operators are intuitively well understood from a semantical perspective via their canonical interpretation over modal Kripke models which formally captures the concept of announcing the truth of a formula to a set of agents by removing the worlds disagreeing with the announcement or appropriately restricting the accessibility relation.

Many-valued versions of public announcement logic, over a t-norm based logic, have only been considered in \cite{CRR2016} where the authors study n-valued {\L}ukasiewicz logic as an underlying base logic and mostly place emphasis on algebraic semantics for public announcements and duality theory. Their relational semantics is, modulo some technicalities, however similar to the \emph{announcement of full truth} considered in Section \ref{sec:fulltruth} later.

The base logic chosen here is G\"odel logic, originating from a note of G\"odel \cite{Goe1932} where he introduced a family of finite-valued logics to show that intuitionistic logic does not have a finite characteristic matrix. These were later extended to an infinitely-valued variant by Dummett \cite{Dum1959} which has been extensively studied and in particular also has been identified as one of the three main t-norm based fuzzy logics by Haj\'ek \cite{Haj1998}.

In particular, G\"odel logic has well-behaved modal extensions as considered in \cite{CR2010} which allow for complete proof calculi with respect to a very natural semantics, generalizing the approach to classical modal logics via Kripke models to truth-values in $[0,1]$. In that context, one usually separates versions with $\Box$-style, $\Diamond$-style and both modalities (see \cite{CR2015} for the latter). We only focus on the $\Box$-fragment for this note.

We consider three possible semantic approaches for a public announcement operator in this context of the $\Box$-fragment of G\"odel modal logic, based on generalizing certain properties of the semantic interpretation of public announcements in the classical setting to the many-valued one. We find that two of those approaches yield different logics in terms of expressive strength and provide indications for this also being true with the third variant. We also provide proof calculi for two of the three introduced semantics.
\subsection{Standard G\"odel modal logics}
We consider the modal language $\mathcal{L}_\Box$ defined via
\[
\mathcal{L}_\Box:\phi::=\bot\mid p\mid (\phi\land\phi)\mid (\phi\rightarrow\phi)\mid\Box\phi.
\]
It shall be noted that we only consider the single modality $\Box$ throughout the paper instead of a family of operators, indexed by a set of agents. However, the results presented here naturally extend to a language with finitely many agents.

Over this language, we consider the calculus of Caicedo and Rodriguez for the $\Box$-fragment of standard G\"odel modal logic from \cite{CR2010}:
\begin{definition}
We define the calculus $\mathcal{GK}$ over the language $\mathcal{L}_\Box$ as follows:
\begin{description}
\item [($G$)] all the axiom schemes of the calculus $\mathcal{G}$;
\item [($K$)] $\Box(\phi\rightarrow\psi)\rightarrow (\Box\phi\rightarrow\Box\psi)$;
\item [($Z$)] $\neg\neg\Box\phi\rightarrow\Box\neg\neg\phi$;
\item [($N\Box$)] from $\vdash\phi$, infer $\vdash\Box\phi$;
\item [($MP$)] from $\phi$ and $\phi\rightarrow\psi$, infer $\psi$.
\end{description}
\end{definition}
Note that $\mathcal{GK}$ has the classical deduction theorem:
\begin{lemma}[\cite{CR2010}]
For any $\Gamma\cup\{\phi,\psi\}\subseteq\mathcal{L}_\Box$: $\Gamma\cup\{\phi\}\vdash_\mathcal{GK}\psi$ iff $\Gamma\vdash_\mathcal{GK}\phi\rightarrow\psi$.
\end{lemma}
Semantically, this logic is characterized by so-called G\"odel-Kripke models, algebraic generalizations of Kripke models to take values in $[0,1]$. In that context, we consider $[0,1]$ as a Heyting algebra with the following functions taking the role of meet and its residuum: we write $x\odot y:=\min\{x,y\}$ for the minimum t-norm for $x,y\in [0,1]$ and $\Rightarrow$ for its residuum, that is
\[
x\Rightarrow y:=\begin{cases}1&\text{if }x\leq y,\\y&\text{if }x>y,\end{cases}
\]
for $x,y\in [0,1]$.
\begin{definition}[\cite{CR2010}]
A G\"odel-Kripke model is a triple $\mathfrak{M}=\langle W,R,V\rangle$ with
\begin{enumerate}
\item $W$ is a non-empty set (the domain $\mathcal{D}(\mathfrak{M})$),
\item $R:W\times W\to [0,1]$,
\item $V:W\times Var\to [0,1]$.
\end{enumerate}
\end{definition}
The function $V(w,\cdot)$ naturally extends, for every $w\in W$, to a function $\vert\cdot\vert^w_\mathfrak{M}:\mathcal{L}_\Box\to [0,1]$ via the following clauses:
\begin{itemize}
\item $\vert\bot\vert^w_\mathfrak{M}:=0$;
\item $\vert p\vert^w_\mathfrak{M}:=V(w,p)$;
\item $\vert\phi\land\psi\vert^w_\mathfrak{M}:=\vert\phi\vert^w_\mathfrak{M}\odot\vert\psi\vert^w_\mathfrak{M}$;
\item $\vert\phi\rightarrow\psi\vert^w_\mathfrak{M}:=\vert\phi\vert^w_\mathfrak{M}\Rightarrow\vert\psi\vert^w_\mathfrak{M}$;
\item $\vert\Box\phi\vert^w_\mathfrak{M}:=\inf\{R(w,v)\Rightarrow\vert\phi\vert^v_\mathfrak{M}\mid v\in W\}$.
\end{itemize}
We extend $\vert\cdot\vert^w_\mathfrak{M}$ to sets of formulas $\Gamma$ via $\vert\Gamma\vert^w_\mathfrak{M}:=\inf\{\vert\gamma\vert^w_\mathfrak{M}\mid \gamma\in\Gamma\}$ and write $(\mathfrak{M},w)\models\phi$ if $\vert\phi\vert^w_\mathfrak{M}=1$ and similarly for sets.

Further, we denote the class of all G\"odel-Kripke models by $\mathsf{GK}$. Over $\mathsf{GK}$, there are now two natural notions of consequence as defined in \cite{CR2010}.
\begin{definition}
Let $\Gamma\cup\{\phi\}\subseteq\mathcal{L}_\Box$. Then, we write
\begin{enumerate}
\item $\Gamma\models_\mathsf{GK}\phi$ if $\forall\mathfrak{M}\in\mathsf{GK}\forall w\in\mathcal{D}(\mathfrak{M})\left((\mathfrak{M},w)\models\Gamma\text{ implies }(\mathfrak{M},w)\models\phi\right)$,
\item $\Gamma\models^\leq_\mathsf{GK}\phi$ if $\forall\mathfrak{M}\in\mathsf{GK}\forall w\in\mathcal{D}(\mathfrak{M})\left(\vert\Gamma\vert^w_\mathfrak{M}\leq\vert\phi\vert^w_\mathfrak{M}\right)$.
\end{enumerate}
\end{definition}
Over those consequence relations, Caicedo and Rodriguez then obtained the following completeness theorem.
\begin{theorem}[\cite{CR2010}]\label{thm:gkcomp}
For any $\Gamma\cup\{\phi\}\subseteq\mathcal{L}_\Box$, the following are equivalent:
\begin{enumerate}
\item $\Gamma\vdash_{\mathcal{GK}}\phi$;
\item $\Gamma\models^\leq_\mathsf{GK}\phi$;
\item $\Gamma\models_\mathsf{GK}\phi$.
\end{enumerate}
\end{theorem}
\subsection{Standard G\"odel modal logic with $\Delta$ and rational constants}\label{sec:SGMLwithDelta}
A common extension of G\"odel (modal) logic is obtained by adding the Baaz $\Delta$-operator \cite{Baa1996} and rational constants to the language. Precisely, we consider
\[
\mathcal{L}^\Delta_\Box([0,1]_\mathbb{Q}):\phi::=\bot\mid\bar{c}\mid p\mid (\phi\land\phi)\mid (\phi\rightarrow\phi)\mid\Delta\phi\mid\Box\phi.
\]
We write $\mathcal{L}^\Delta_\Box$ for the sublanguage without rational constants. 

Semantically, the same G\"odel-Kripke models apply here, where the evaluation functions are extended by the clauses
\[
\vert	\bar{c}\vert^w_\mathfrak{M}:=c\text{ and }\vert\Delta\phi\vert^w_\mathfrak{M}:=\delta(\vert\phi\vert^w_\mathfrak{M})
\]
where $\delta$ is defined by
\[
\delta(x):=\begin{cases}1&\text{if }x=1,\\0&\text{otherwise},\end{cases}
\]
for $x\in [0,1]$. The consequence relation $\models_\mathsf{GK}$ then naturally extends to the new language.

A proof calculus for the theory of G\"odel-Kripke models over that language is obtained by extending $\mathcal{GK}$ with axioms for $\Delta$ and the constants together with some connection axioms and is given in the following definition.
\begin{definition}
We define the calculus $\mathcal{GK}_\Delta([0,1]_\mathbb{Q})$ as follows:
\begin{description}
\item [($GK$)] all the axiom schemes of the calculus $\mathcal{GK}$;
\item [($\Delta 1$)] $\Delta\phi\lor\neg\Delta\phi$;
\item [($\Delta 2$)] $\Delta(\phi\lor\psi)\rightarrow\Delta\phi\lor\Delta\psi$;
\item [($\Delta 3$)] $\Delta\phi\rightarrow\phi$;
\item [($\Delta 4$)] $\Delta\phi\rightarrow\Delta\Delta\phi$;
\item [($\Delta 5$)] $\Delta(\phi\rightarrow\psi)\rightarrow(\Delta\phi\rightarrow\Delta\psi)$;
\item [($BK1$)] $\bar{c}\land\bar{d}\leftrightarrow\overline{c\odot d}$;
\item [($BK2$)] $\bar{c}\rightarrow\bar{d}\leftrightarrow\overline{c\Rightarrow d}$;
\item [($BK3$)] $\Delta\bar{c}\leftrightarrow\bar{\delta(c)}$;
\item [($BK4$)] $\bot\leftrightarrow\bar{0}$;
\item [($\Delta\Box$)] $\Delta\Box\phi\rightarrow\Box\Delta\phi$;
\item [($\Box\bar{c}$)] $\Box(\bar{c}\rightarrow\phi)\leftrightarrow\bar{c}\rightarrow\Box\phi$;
\item [($N\Delta$)] from $\phi$, infer $\Delta\phi$;
\item [($N\Box$)] from $\vdash\phi$, infer $\vdash\Box\phi$;
\item [($MP$)] from $\phi\rightarrow\psi$ and $\phi$, infer $\psi$.
\end{description}
\end{definition}
Indeed, this calculus is strongly complete w.r.t. the semantics over G\"odel-Kripke models as shown by Vidal in \cite{Vid2015} by reducing it to the underlying propositional G\"odel logic with rational constants and $\Delta$ (see \cite{EGGN2007} for a comprehensive study of such extensions in the context of t-norm based fuzzy logics):
\begin{theorem}[\cite{Vid2015}]
For any $\Gamma\cup\{\phi\}\subseteq\mathcal{L}_\Box^\Delta([0,1]_\mathbb{Q})$, the following are equivalent:
\begin{enumerate}
\item $\Gamma\vdash_{\mathcal{GK}_\Delta([0,1]_\mathbb{Q})}\phi$,
\item $\Gamma\models^\leq_\mathsf{GK}\phi$,
\item $\Gamma\models_\mathsf{GK}\phi$.
\end{enumerate}
\end{theorem}
\section{The good: a first public announcement operator}\label{sec:simplePAop}
We begin by extending the language $\mathcal{L}_\Box$ to a language $\mathcal{L}_{PA}$ by adding a usual syntactic public announcement operator $[\cdot]$:
\[
\mathcal{L}_{PA}:\phi::=\bot\mid p\mid (\phi\land\phi)\mid (\phi\rightarrow\phi)\mid\Box\phi\mid [\phi]\phi.
\]
In this section, we will study a possible semantics for a naive generalization of the classical public announcement logic over $\mathcal{K}$. This \emph{naive} generalization arises by taking a usual axiomatization of classical public announcement logic and replacing the base of the classical modal logic $\mathcal{K}$ by the above $\mathcal{GK}_\Box$. In other words, as classical public announcement logic is axiomatized by reduction axioms for the interplay between the public announcement operator and the other connectives, the calculus we study looks like this.
\begin{definition}
We define the calculus $\mathcal{GPA}$ over the language $\mathcal{L}_\Box$ as follows:
\begin{description}
\item [($GK$)] all the axiom schemes and rules of the calculus $\mathcal{GK}$;
\item [($PA1$)] $[\phi]\bot\leftrightarrow(\phi\rightarrow\bot)$;
\item [($PA2$)] $[\phi]p\leftrightarrow(\phi\rightarrow p)$;
\item [($PA3$)] $[\phi](\psi\land\chi)\leftrightarrow([\phi]\psi\land[\phi]\chi)$;
\item [($PA4$)] $[\phi](\psi\rightarrow\chi)\leftrightarrow([\phi]\psi\rightarrow[\phi]\chi)$;
\item [($PA5$)] $[\phi]\Box\psi\leftrightarrow(\phi\rightarrow\Box[\phi]\psi)$;
\item [($PA6$)] $[\phi][\psi]\chi\leftrightarrow[\phi\land[\phi]\psi]\chi$.
\end{description}
\end{definition}
\subsection{Semantics}
The corresponding semantics is then obtained by considering the same G\"odel-Kripke models as before, where the evaluation function is extended to $\mathcal{L}_{PA}$ by adding the clause
\[
\vert [\phi]\psi\vert^w_\mathfrak{M}:=\vert\phi\vert^w_\mathfrak{M}\Rightarrow\vert\psi\vert^w_{\mathfrak{M}\vert\phi}
\]
for a G\"odel-Kripke model $\mathfrak{M}=\langle W,R,V\rangle$ with $w\in W$ and where the model $\mathfrak{M}\vert\phi:=\langle W^\phi_{\mathfrak{M}},R^\phi_{\mathfrak{M}},V^\phi_{\mathfrak{M}}\rangle$ is defined by $W^\phi_{\mathfrak{M}}:= W$ and $V^\phi_{\mathfrak{M}}:=V$ as well as
\[
R^\phi_{\mathfrak{M}}(w,v):=R(w,v)\odot\vert\phi\vert^v_\mathfrak{M}.
\]
Via this extended evaluation function, the consequence relation $\models_\mathsf{GK}$ naturally extends to inputs $\Gamma\cup\{\phi\}\subseteq\mathcal{L}_{PA}$.
\begin{lemma}\label{lem:fpavalidity}
The following formulas are valid in all G\"odel-Kripke models:
\begin{enumerate}
\item $[\phi]\bot\leftrightarrow(\phi\rightarrow\bot)$;
\item $[\phi]p\leftrightarrow(\phi\rightarrow p)$;
\item $[\phi](\psi\land\chi)\leftrightarrow([\phi]\psi\land[\phi]\chi)$;
\item $[\phi](\psi\rightarrow\chi)\leftrightarrow([\phi]\psi\rightarrow[\phi]\chi)$;
\item $[\phi]\Box\psi\leftrightarrow(\phi\rightarrow\Box[\phi]\psi)$;
\item $[\phi][\psi]\chi\leftrightarrow[\phi\land[\phi]\psi]\chi$.
\end{enumerate}
\end{lemma}
\begin{proof}
In the following, let $\mathfrak{M}=\langle W,R,V\rangle$ be a G\"odel-Kripke model and let $w\in W$. Items (1) and (2) are rather immediate.
\begin{enumerate}
\setcounter{enumi}{2}
\item We have
\begin{align*}
\vert[\phi](\psi\land\chi)\vert^w_\mathfrak{M}
&=\vert\phi\vert^w_\mathfrak{M}\Rightarrow\left(\vert\psi\vert^w_{\mathfrak{M}\vert\phi}\odot\vert\chi\vert^w_{\mathfrak{M}\vert\phi}\right)\\
&=\left(\vert\phi\vert^w_\mathfrak{M}\Rightarrow\vert\psi\vert^w_{\mathfrak{M}\vert\phi}\right)\odot\left(\vert\phi\vert^w_\mathfrak{M}\Rightarrow\vert\chi\vert^w_{\mathfrak{M}\vert\phi}\right)\\
&=\vert[\phi]\psi\vert^w_\mathfrak{M}\odot\vert[\phi]\chi\vert^w_\mathfrak{M}
\end{align*}
where we use the $\mathbf{[0,1]_G}$-identity (even valid in general Heyting algebras)
\[
x\Rightarrow (y\odot z)=(x\Rightarrow y)\odot (x\Rightarrow z)
\]
for all $x,y,z\in [0,1]$.
\item We have
\begin{align*}
\vert[\phi](\psi\rightarrow\chi)\vert^w_\mathfrak{M}
&=\vert\phi\vert^w_\mathfrak{M}\Rightarrow\left(\vert\psi\vert^w_{\mathfrak{M}\vert\phi}\Rightarrow\vert\chi\vert^w_{\mathfrak{M}\vert\phi}\right)\\
&=\left(\vert\phi\vert^w_\mathfrak{M}\Rightarrow\vert\psi\vert^w_{\mathfrak{M}\vert\phi}\right)\Rightarrow\left(\vert\phi\vert^w_\mathfrak{M}\Rightarrow\vert\chi\vert^w_{\mathfrak{M}\vert\phi}\right)\\
&=\vert[\phi]\psi\vert^w_\mathfrak{M}\Rightarrow\vert[\phi]\chi\vert^w_\mathfrak{M}
\end{align*}
where we use the $\mathbf{[0,1]_G}$-identity
\[
x\Rightarrow (y\Rightarrow z)=(x\Rightarrow y)\Rightarrow (x\Rightarrow z)
\]
for all $x,y,z\in [0,1]$.
\item We have
\begin{align*}
\vert[\phi]\Box\psi\vert^w_\mathfrak{M}
&=\vert\phi\vert^w_\mathfrak{M}\Rightarrow\left(\vert\Box\psi\vert^w_{\mathfrak{M}\vert\phi}\right)\\
&=\vert\phi\vert^w_\mathfrak{M}\Rightarrow\left(\inf_{v\in W^\phi_\mathfrak{M}}\left\{R^\phi_\mathfrak{M}(w,v)\Rightarrow\vert\psi\vert^v_{\mathfrak{M}\vert\phi}\right\}\right)\\
&=\vert\phi\vert^w_\mathfrak{M}\Rightarrow\left(\inf_{v\in W}\left\{\left(R(w,v)\odot\vert\phi\vert^v_\mathfrak{M}\right)\Rightarrow\vert\psi\vert^v_{\mathfrak{M}\vert\phi}\right\}\right)\\
&=\vert\phi\vert^w_\mathfrak{M}\Rightarrow\left(\inf_{v\in W}\left\{R(w,v)\Rightarrow\left(\vert\phi\vert^v_\mathfrak{M}\Rightarrow\vert\psi\vert^v_{\mathfrak{M}\vert\phi}\right)\right\}\right)\\
&=\vert\phi\vert^w_\mathfrak{M}\Rightarrow\left(\inf_{v\in W}\left\{R(w,v)\Rightarrow\vert[\phi]\psi\vert^v_\mathfrak{M}\right\}\right)\\
&=\vert\phi\vert^w_\mathfrak{M}\Rightarrow\vert\Box[\phi]\psi\vert^w_\mathfrak{M}
\end{align*}
where, at the fourth equality, we used the $\mathbf{[0,1]_G}$-identity
\[
(x\odot y)\Rightarrow z=x\Rightarrow (y\Rightarrow z)
\]
for $x,y,z\in [0,1]$.
\item We have
\begin{align*}
\vert[\phi][\psi]\chi\vert^w_\mathfrak{M}
&=\vert\phi\vert^w_\mathfrak{M}\Rightarrow\vert [\psi]\chi\vert^w_{\mathfrak{M}\vert\phi}\\
&=\vert\phi\vert^w_\mathfrak{M}\Rightarrow\left(\vert\psi\vert^w_{\mathfrak{M}\vert\phi}\Rightarrow \vert\chi\vert^w_{(\mathfrak{M}\vert\phi)\vert\psi}\right)\\
&=\left(\vert\phi\vert^w_\mathfrak{M}\odot\vert\psi\vert^w_{\mathfrak{M}\vert\phi}\right)\Rightarrow \vert\chi\vert^w_{(\mathfrak{M}\vert\phi)\vert\psi}
\end{align*}
using the previous identity 
\[
(x\odot y)\Rightarrow z=x\Rightarrow (y\Rightarrow z).
\]
Also, we get
\begin{align*}
\vert[\phi\land[\phi]\psi]\chi\vert^w_\mathfrak{M}
&=\vert\phi\land[\phi]\psi\vert^w_\mathfrak{M}\Rightarrow\vert\chi\vert^w_{\mathfrak{M}\vert(\phi\land[\phi]\psi)}\\
&=\left(\vert\phi\vert^w_\mathfrak{M}\odot\vert[\phi]\psi\vert^w_\mathfrak{M}\right)\Rightarrow\vert\chi\vert^w_{\mathfrak{M}\vert(\phi\land[\phi]\psi)}\\
&=\left(\vert\phi\vert^w_\mathfrak{M}\odot\left(\vert\phi\vert^w_\mathfrak{M}\Rightarrow \vert\psi\vert^w_{\mathfrak{M}\vert\phi}\right)\right)\Rightarrow\vert\chi\vert^w_{\mathfrak{M}\vert(\phi\land[\phi]\psi)}\\
&=\left(\vert\phi\vert^w_\mathfrak{M}\odot\vert\psi\vert^w_{\mathfrak{M}\vert\phi}\right)\Rightarrow\vert\chi\vert^w_{\mathfrak{M}\vert(\phi\land[\phi]\psi)}
\end{align*}
were now, the contributing identity of $\mathbf{[0,1]_G}$ is
\[
x\odot (x\Rightarrow y)=x\odot y
\]
for $x,y\in [0,1]$.

So, the two parts are equal if
\[
\vert\chi\vert^w_{\mathfrak{M}\vert(\phi\land[\phi]\psi)}=\vert\chi\vert^w_{(\mathfrak{M}\vert\phi)\vert\psi}.
\]
But this is satisfied as 
\[
\mathfrak{M}\vert(\phi\land[\phi]\psi)=(\mathfrak{M}\vert\phi)\vert\psi.
\]
For that, it suffices to note that
\begin{align*}
R_\mathfrak{M}^{(\phi\land[\phi]\psi)}(w,v)
&=R(w,v)\odot\vert\phi\land[\phi]\psi\vert^v_\mathfrak{M}\\
&=R(w,v)\odot\left(\vert\phi\vert^v_\mathfrak{M}\odot\left(\vert\phi\vert^v_\mathfrak{M}\Rightarrow\vert\psi\vert^v_{\mathfrak{M}\vert\phi}\right)\right)\\
&=R(w,v)\odot\vert\phi\vert^v_\mathfrak{M}\odot\vert\psi\vert^v_{\mathfrak{M}\vert\phi}\\
&=R_\mathfrak{M}^\phi(w,v)\odot\vert\psi\vert^v_{\mathfrak{M}\vert\phi}\\
&=\left(R_\mathfrak{M}^\phi\right)_{\mathfrak{M}\vert\phi}^\psi(w,v).
\end{align*}
\end{enumerate}
\end{proof}
From this, the following lemma is rather immediate.
\begin{lemma}\label{lem:pasoundness}
For any $\Gamma\cup\{\phi\}\subseteq\mathcal{L}_{PA}$: $\Gamma\vdash_{\mathcal{GPA}}\phi$ implies $\Gamma\models^\leq_\mathsf{GK}\phi$.
\end{lemma}
Now, the following completeness proof is standard as we reduce, in the same way as in the classical case, formulas containing public announcement operators to equivalent formulas from $\mathcal{L}_\Box$ alone.
\begin{definition}
We define the function $t:\mathcal{L}_{PA}\to\mathcal{L}_\Box$ by recursion as follows:
\begin{enumerate}
\item $t(p)=p$ for $p\in Var$; $t(\bot)=\bot$;
\item $t(\phi\circ\psi)=t(\phi)\circ t(\psi)$ for $\circ\in\{\land,\rightarrow\}$;
\item $t(\Box\phi)=\Box t(\phi)$;
\item $t([\phi]\bot)=t(\phi)\rightarrow\bot$;
\item $t([\phi]p)=t(\phi)\rightarrow p$;
\item $t([\phi](\psi\land\chi))=t([\phi]\psi)\land t([\phi]\chi)$;
\item $t([\phi](\psi\rightarrow\chi))=t([\phi]\psi)\rightarrow t([\phi]\chi)$;
\item $t([\phi]\Box\psi)=t(\phi)\rightarrow\Box t([\phi]\psi)$;
\item $t([\phi][\psi]\chi)=t([\phi\land[\phi]\psi]\chi)$.
\end{enumerate}
\end{definition}
This translation provides equivalent formulas which contain no public announcement and this equivalence is provable in $\mathcal{GPA}$:
\begin{lemma}\label{lem:transval}
For all $\phi\in\mathcal{L}_{PA}$: $\vdash_{\mathcal{GPA}}\phi\leftrightarrow t(\phi)$.
\end{lemma}
The proof follows the standard procedure of defining a suitable complexity measure over which we then perform an induction (see e.g. \cite{DHK2008} for that in the classical case). We include a sketch here for self-containedness.
\begin{proof}
As said in the paragraph above, the key element is a complexity measure on $\mathcal{L}_{PA}$ such that one can perform a suitable induction. Following \cite{DHK2008}, we use the function $c:\mathcal{L}_{PA}\to\mathbb{N}$ defined via
\begin{itemize}
\item $c(p)=c(\bot)=1$,
\item $c(\phi\circ\psi)=1+\max\{c(\phi),c(\psi)\}$ for $\circ\in\{\land,\rightarrow\}$,
\item $c(\Box\phi)=1+c(\phi)$,
\item $c([\phi]\psi)=(4+c(\phi))\cdot c(\psi)$.
\end{itemize}
This $c$ can be easily seen to have the following properties (similar as in \cite{DHK2008}):
\begin{enumerate}
\item $c(\psi)\leq c(\phi)$ for $\psi$ a subformula of $\phi$;
\item $c(\phi\rightarrow\bot)<c([\phi]\bot)$;
\item $c(\phi\rightarrow p)<c([\phi]p)$;
\item $c([\phi]\psi\land[\phi]\chi)<c([\phi](\psi\land\chi))$;
\item $c([\phi]\psi\rightarrow[\phi]\chi)<c([\phi](\psi\rightarrow\chi))$;
\item $c(\phi\rightarrow\Box[\phi]\psi)<c([\phi]\Box\psi)$;
\item $c([\phi\land[\phi]\psi]\chi)<c([\phi][\psi]\chi)$.
\end{enumerate}
Using the reduction axioms, it is now straightforward to prove the theorem by induction on $c(\phi)$.
\end{proof}
This translation, together with soundness, implies the expressive equivalence between the standard G\"odel modal logic and its extension by public announcements.
\begin{theorem}
$\mathcal{L}_\Box\equiv\mathcal{L}_{PA}$ over all G\"odel-Kripke models.
\end{theorem}
Further, we can actually use this translation, as coded into the proof calculus $\mathcal{GPA}$ by the reduction axioms, to provide a completeness result.
\begin{theorem}\label{thm:pacomp}
For any $\Gamma\cup\{\phi\}\subseteq\mathcal{L}_{PA}$, the following are equivalent:
\begin{enumerate}
\item $\Gamma\vdash_{\mathcal{GPA}}\phi$;
\item $\Gamma\models^\leq_{\mathsf{GK}}\phi$;
\item $\Gamma\models_{\mathsf{GK}}\phi$.
\end{enumerate}
\end{theorem}
\begin{proof}
``(1) implies (2)" follows from Lemma \ref{lem:pasoundness} and ``(2) implies (3)" is immediate. For ``(3) implies (1)", suppose $\Gamma\models_\mathsf{GK}\phi$. By Lemma \ref{lem:transval} together with Lemma \ref{lem:pasoundness}, we have
\[
\models_\mathsf{GK}\psi\leftrightarrow t(\psi)
\]
for all $\psi\in\mathcal{L}_{PA}$ and therefore 
\[
t[\Gamma]\models_\mathsf{GK}t(\phi).
\]
Completeness of $\mathsf{GK}$, Theorem \ref{thm:gkcomp}, implies $t[\Gamma]\vdash_\mathcal{GK}t(\phi)$ and thus
\[
t[\Gamma]\vdash_\mathcal{GPA}t(\phi).
\]
Now, Lemma \ref{lem:transval} immediately yields $\Gamma\vdash_\mathcal{GPA}\phi$.
\end{proof}
\section{The bad: announcing full truth}\label{sec:fulltruth}
Another possible fuzzy public announcement operator is obtained by announcing that $\phi$ has truth value $1$. To distinguish this from the previous operator, we consider the language
\[
\mathcal{L}_{PA^1}:\phi::=\bot\mid p\mid (\phi\land\phi)\mid (\phi\rightarrow\phi)\mid\Box\phi\mid [\phi]^1\phi.
\]
Semantically, the language is again interpreted over G\"odel-Kripke model via the additional clause 
\[
\vert [\phi]^1\psi\vert^w_\mathfrak{M}:=\begin{cases}\vert\psi\vert^w_{\mathfrak{M}\vert^1\phi}&\text{if }\vert\phi\vert^w_\mathfrak{M}=1,\\
1&\text{otherwise},
\end{cases}
\]
for a G\"odel-Kripke model $\mathfrak{M}=\langle W,R,V\rangle$ with $w\in W$, but where the model $\mathfrak{M}\vert^1\phi:=\langle W_\mathfrak{M}^{\phi,1},R_\mathfrak{M}^{\phi,1},V_\mathfrak{M}^{\phi,1}\rangle$ is now defined by restricting the set of worlds $W$ to
\[
W_\mathfrak{M}^{\phi,1}:=\{v\in W\mid\vert\phi\vert^v_\mathfrak{M}=1\}
\] 
and $V_\mathfrak{M}^{\phi,1}$ as well as $R_\mathfrak{M}^{\phi,1}$ are just induced from $V$ and $R$ by this new set of worlds. Again, this extended evaluation function allows us to lift the consequence relation $\models_\mathsf{GK}$ to inputs $\Gamma\cup\{\phi\}\subseteq\mathcal{L}_{PA^1}$.\\

Now, this semantics turns out to be intimately connected with the $\Delta$-operator from \cite{Baa1996} (see also Section \ref{sec:SGMLwithDelta} again). This is preliminary exemplified by the following lemma which present modified reduction axioms.\footnote{The $\Delta$-operator is not officially part of the language $\mathcal{L}_{PA^1}$ but we pretend so for convenience in the following lemma.}
\begin{lemma}\label{lem:1paaxiomsvalid}
The following schemes are valid in any $\mathsf{GK}$-model:
\begin{enumerate}
\item $[\phi]^1\bot\leftrightarrow(\Delta\phi\rightarrow\bot)$;
\item $[\phi]^1p\leftrightarrow(\Delta\phi\rightarrow p)$;
\item $[\phi]^1(\psi\land\chi)\leftrightarrow([\phi]^1\psi\land[\phi]^1\chi)$;
\item $[\phi]^1(\psi\rightarrow\chi)\leftrightarrow([\phi]^1\psi\rightarrow[\phi]^1\chi)$;
\item $[\phi]^1\Box\psi\leftrightarrow(\Delta\phi\rightarrow\Box[\phi]^1\psi)$;
\item $[\phi]^1[\psi]^1\chi\leftrightarrow[\phi\land[\phi]^1\psi]^1\chi$.
\end{enumerate}
\end{lemma}
\begin{proof}
The proof is similar to that of Lemma \ref{lem:fpavalidity} and we only show the interesting cases where there are new concepts coming in. In particular, we only show items (1), (5) and (6). For that, let $\mathfrak{M}$ be a G\"odel-Kripke model and let $w\in\mathcal{D}(\mathfrak{M})$.
\begin{enumerate}
\item We have
\begin{align*}
\vert[\phi]^1\bot\vert^w_\mathfrak{M}&=\begin{cases}0&\text{if }\vert\phi\vert^w_\mathfrak{M}=1\\1&\text{otherwise}\end{cases}\\
&=\sim\delta(\vert\phi\vert^w_\mathfrak{M}).
\end{align*}
\setcounter{enumi}{4}
\item We have
\begin{align*}
\vert[\phi]^1\Box\psi\vert^w_\mathfrak{M}
&=\begin{cases}\vert\Box\psi\vert^w_{\mathfrak{M}\vert^1\phi}&\text{if }\vert\phi\vert^w_\mathfrak{M}=1\\
1&\text{otherwise}\end{cases}\\
&=\begin{cases}\inf_{v\in W, \vert\phi\vert^v_\mathfrak{M}=1}\left\{R(w,v)\Rightarrow\vert\psi\vert^v_{\mathfrak{M}\vert^1\phi}\right\}&\text{if }\vert\phi\vert^w_\mathfrak{M}=1\\
1&\text{otherwise}\end{cases}\\
&=\vert\Delta\phi\vert^w_\mathfrak{M}\Rightarrow\inf_{v\in W, \vert\phi\vert^v_\mathfrak{M}=1}\left\{R(w,v)\Rightarrow\vert\psi\vert^v_{\mathfrak{M}\vert^1\phi}\right\}\\
&=\vert\Delta\phi\vert^w_\mathfrak{M}\Rightarrow\left(\inf_{v\in W, \vert\phi\vert^v_\mathfrak{M}=1}\left\{R(w,v)\Rightarrow\vert\psi\vert^v_{\mathfrak{M}\vert^1\phi}\right\}\odot \inf_{v\in W, \vert\phi\vert^v_\mathfrak{M}\neq 1}\left\{R(w,v)\Rightarrow 1\right\}\right)\\
&=\vert\Delta\phi\vert^w_\mathfrak{M}\Rightarrow\inf_{v\in W}\left\{R(w,v)\Rightarrow\vert[\phi]^1\psi\vert^v_{\mathfrak{M}}\right\}\\
&=\vert\Delta\phi\rightarrow\Box[\phi]^1\psi\vert^v_\mathfrak{M}.
\end{align*}
\item We first show $(W\vert^1\phi)\vert^1\psi=W\vert^1(\phi\land[\phi]^1\psi)$. To see that, note 
\begin{align*}
(W\vert^1\phi)\vert^1\psi&=\{w\in W\vert^{\phi,1}_\mathfrak{M}\mid\vert\psi\vert^w_{\mathfrak{M}\vert^1\phi}=1\}\\
&=\{w\in W\mid \vert\phi\vert^w_\mathfrak{M}=1\text{ and }\vert\psi\vert^w_{\mathfrak{M}\vert^1\phi}=1\}\\
&=\{w\in W\mid \vert\phi\vert^w_\mathfrak{M}=1\text{ and }\vert[\phi]^1\psi\vert^w_\mathfrak{M}=1\}\\
&=\{w\in W\mid \vert\phi\land[\phi]^1\psi\vert^w_\mathfrak{M}=1\}\\
&=W\vert^1(\phi\land[\phi]^1\psi).
\end{align*}
Now, using the equality, we can see the validity of item (8) rather quickly:
\begin{align*}
\vert[\phi]^1[\psi]^1\chi\vert^w_\mathfrak{M}&=\begin{cases}\vert\chi\vert^w_{(\mathfrak{M}\vert^1\phi)\vert^1\psi}&\text{if }\vert\phi\vert^w_\mathfrak{M}=1\text{ and }\vert\psi\vert^w_{\mathfrak{M}\vert^1\phi}=1,\\1&\text{otherwise},\end{cases}\\
&=\begin{cases}\vert\chi\vert^w_{\mathfrak{M}\vert^1(\phi\land[\phi]\psi)}&\text{if }\vert\phi\vert^w_\mathfrak{M}=1\text{ and }\vert[\phi]^1\psi\vert^w_{\mathfrak{M}}=1,\\1&\text{otherwise},\end{cases}\\
&=\vert[\phi\land[\phi]^1\psi]^1\chi\vert^w_\mathfrak{M}.
\end{align*}
\end{enumerate}
\end{proof}
But further, we can actually define the $\Delta$ operator via 1-announcements. This is being hinted in item (1) of the previous lemma which implies that
\[
\neg\Delta\phi\leftrightarrow[\phi]^1\bot
\]
is valid in all $\mathsf{GK}$-models. Now, normally $\neg\neg\phi$ is not equivalent to $\phi$ in G\"odel (modal) logic but indeed this holds for $\Delta\psi$. So, we obtain
\[
\Delta\phi\leftrightarrow\neg\neg\Delta\phi\leftrightarrow\neg [\phi]^1\bot.
\]
Now, this observation is collected in the following lemma.
\begin{lemma}
The scheme $\neg[\phi]^1\bot\leftrightarrow\Delta\phi$ is valid in any $\mathsf{GK}$-model.
\end{lemma}
The above results immediately yield translations between $\mathcal{L}^\Delta_\Box$ and $\mathcal{L}_{PA^1}$ (and $\mathcal{L}^\Delta_{PA^1}$, that is $\mathcal{L}_{PA^1}$ extended by $\Delta$) which provide equal expressivity.
\begin{theorem}
$\mathcal{L}_{PA^1}\equiv\mathcal{L}^\Delta_\Box\equiv\mathcal{L}^\Delta_{PA^1}$.
\end{theorem}
\subsection{Adding rational constants and the $\Delta$-operator}
Now, the outline above, and in particular the expressivity result, may make it compelling to provide a proof calculus first for $\mathcal{L}^\Delta_{PA^1}$ by using the definable translation into $\mathcal{L}^\Delta_\Box$ from the previous lemmas to provide reduction axioms as in the case of $\mathcal{L}_{PA}$ and $\mathcal{L}_\Box$. Then, one may actually additionally drop the $\Delta$ and obtain an axiomatization for $\mathcal{L}_{PA^1}$ by using the internal definability of $\Delta$ within $\mathcal{L}_{PA^1}$.

We are however not aware of any complete proof calculus for the semantic consequence of G\"odel-Kripke models over $\mathcal{L}^\Delta_\Box$ and an inspection of the usual completeness proof of Caicedo and Rodriguez for basic G\"odel modal logic as given in \cite{CR2010} makes it apparent that their approach does not work in the context of the $\Delta$-operator as the underlying propositional G\"odel logic with $\Delta$ does not have the classical deduction theorem. If one could provide such a proof calculus for $\mathcal{L}^\Delta_\Box$, the above would immediately yield, as outlined, proof calculi for both $\mathcal{L}^\Delta_{PA^1}$ and $\mathcal{L}_{PA^1}$.

There is, however, a marginal extension of $\mathcal{L}^\Delta_\Box$ such that the resulting semantics over G\"odel-Kripke models can be completely axiomatized. More precisely, we may add explicit constants for rational values in $[0,1]$ into the language as has been carried out by Vidal \cite{Vid2015}. This semantics has already been presented in Section \ref{sec:SGMLwithDelta}.

We thus consider the extension of the above language by $\Delta$ and by rational constants:
\[
\mathcal{L}^\Delta_{PA^1}([0,1]_\mathbb{Q}):\phi::=\bot\mid\bar{c}\mid p\mid (\phi\land\phi)\mid (\phi\rightarrow\phi)\mid\Delta\phi\mid\Box\phi\mid [\phi]^1\phi.
\]
This extension indeed does not provide any complications with public announcements, modalities, or the $\Delta$-operator as the following lemma shows.
\begin{lemma}\label{lem:1paaxiomsvalid2}
The following schemes are valid in any $\mathsf{GK}$-model:
\begin{enumerate}
\item $[\phi]^1\bar{c}\leftrightarrow(\Delta\phi\rightarrow \bar{c})$;
\item $[\phi]^1\Delta\psi\leftrightarrow(\Delta\phi\rightarrow\Delta[\phi]^1\psi)$.
\end{enumerate}
\end{lemma}
\begin{proof}
The first item can be shown in the same way as item (1) of Lemma \ref{lem:1paaxiomsvalid}. Regarding the second item, note that 
\[
\vert[\phi]^1\Delta\psi\vert^w_\mathfrak{M}=\begin{cases}\delta(\vert\psi\vert^w_{\mathfrak{M}\vert^1\phi})&\text{if }\vert\phi\vert^w_\mathfrak{M}=1,\\1&\text{otherwise}.\end{cases}
\]
At the same time, we have
\begin{align*}
\vert\Delta\phi\rightarrow\Delta[\phi]^1\psi\vert^w_\mathfrak{M}&=\begin{cases}\delta(\vert [\phi]^1\psi\vert^w_{\mathfrak{M}})&\text{if }\vert\phi\vert^w_\mathfrak{M}=1,\\1&\text{otherwise},\end{cases}\\
&=\begin{cases}\delta(\vert\psi\vert^w_{\mathfrak{M}\vert^1\phi})&\text{if }\vert\phi\vert^w_\mathfrak{M}=1,\\1&\text{otherwise},\end{cases}
\end{align*}
since, if $\vert\phi\vert^w_\mathfrak{M}=1$, then $\vert[\phi]^1\psi\vert^w_\mathfrak{M}=\vert\psi\vert^w_{\mathfrak{M}\vert^1\phi}$.
\end{proof}
Before we provide a proof calculus for the resulting logic, we want to provide a small remedy for the rather unsatisfactory state of not having a proof calculus for $\mathcal{L}_{\Box}^\Delta$ over G\"odel-Kripke models and thus moving to an extension by rational constants: these constants allow for rather nice derived public announcement operators which may be useful in the context of modeling fuzzy epistemic scenarios.

For rational $x$, we may define
\begin{enumerate}
\item $[\phi]^{=x}\psi:=[\Delta(\phi\leftrightarrow\bar{x})]^1\psi$,
\item $[\phi]^{\leq x}\psi:=[\Delta(\phi\rightarrow\bar{x})]^1\psi$,
\item $[\phi]^{>x}\psi:=[\neg\Delta(\phi\rightarrow\bar{x})]^1\psi$,
\item $[\phi]^{\geq x}\psi:=[\neg\Delta(\phi\rightarrow\bar{x})\lor\Delta(\phi\leftrightarrow\bar{x})]^1\psi$,
\item $[\phi]^{<x}\psi:=[\Delta(\phi\rightarrow\bar{x})\land\neg\Delta(\phi\leftrightarrow\bar{x})]^1\psi$
\end{enumerate}
The $\Delta$-operator in (1) and (2) is actually superfluous. Indeed, these operators have the following semantical interpretation.
\begin{lemma}
Let $\mathfrak{M}$ be a G\"odel-Kripke model and $w\in\mathcal{D}(\mathfrak{M})$. Let $\lhd\in\{=,\leq,\geq,<,>\}$ and $x\in[0,1]_\mathbb{Q}$. Then
\[
\vert[\phi]^{\lhd x}\psi\vert^w_\mathfrak{M}=\begin{cases}\vert\psi\vert^w_{\mathfrak{M}\vert^w_{\phi,\lhd x}}&\text{if }\vert\phi\vert^w_\mathfrak{M}\lhd x,\\1&\text{otherwise},\end{cases}
\]
where $\mathfrak{M}\vert^w_{\phi,\lhd x}:=\langle W^{\mathfrak{M},x}_{\phi,\lhd x},R^{\mathfrak{M},w}_{\phi,\lhd x},V^{\mathfrak{M},w}_{\phi,\lhd x}\rangle$ where again $R^{\mathfrak{M},w}_{\phi,\lhd x},V^{\mathfrak{M},w}_{\phi,\lhd x}$ are restrictions and
\[
W^{\mathfrak{M},x}_{\phi,\lhd x}:=\{v\in W\mid\vert\phi\vert^v_\mathfrak{M}\lhd x\}.
\]
\end{lemma}
\subsection{An axiomatization of the extension}
We can now provide an axiomatization for the full language $\mathcal{L}^\Delta_{PA^1}([0,1]_\mathbb{Q})$ over G\"odel-Kripke models by adding the suitable reduction axioms as suggested by the translation hinted in the previous lemmas. The rest of the completeness proof then follows the same lines as before.
\begin{definition}
We define the calculus $\mathcal{GPA}^1_\Delta([0,1]_\mathbb{Q})$ over the language $\mathcal{L}^\Delta_{PA^1}([0,1]_\mathbb{Q})$ as follows:
\begin{description}
\item [($GK'$)] all the axiom schemes and rules of the calculus $\mathcal{GK}_\Delta([0,1]_\mathbb{Q})$;
\item [($PA1$)] $[\phi]^1\bot\leftrightarrow(\Delta\phi\rightarrow\bot)$;
\item [($PA2$)] $[\phi]^1p\leftrightarrow(\Delta\phi\rightarrow p)$;
\item [($PA3$)] $[\phi]^1\bar{c}\leftrightarrow(\Delta\phi\rightarrow \bar{c})$;
\item [($PA4$)] $[\phi]^1(\psi\land\chi)\leftrightarrow([\phi]^1\psi\land[\phi]^1\chi)$;
\item [($PA5$)] $[\phi]^1(\psi\rightarrow\chi)\leftrightarrow([\phi]^1\psi\rightarrow[\phi]^1\chi)$;
\item [($PA6$)] $[\phi]^1\Box\psi\leftrightarrow(\Delta\phi\rightarrow\Box[\phi]^1\psi)$;
\item [($PA7$)] $[\phi]^1\Delta\psi\leftrightarrow(\Delta\phi\rightarrow\Delta[\phi]^1\psi)$;
\item [($PA8$)] $[\phi]^1[\psi]^1\chi\leftrightarrow[\phi\land[\phi]^1\psi]^1\chi$.
\end{description}
\end{definition}
Again, soundness is immediate using Lemma \ref{lem:1paaxiomsvalid} and \ref{lem:1paaxiomsvalid2} and an induction on the length of the proof.
\begin{lemma}\label{lem:1pasoundness}
For any $\Gamma\cup\{\phi\}\subseteq\mathcal{L}^\Delta_{PA^1}([0,1]_\mathbb{Q})$: $\Gamma\vdash_{\mathcal{GPA}^1_\Delta([0,1]_\mathbb{Q})}\phi$ implies $\Gamma\models^\leq_\mathsf{GK}\phi$.
\end{lemma}

\begin{definition}
We define the function $t':\mathcal{L}^\Delta_{PA^1}([0,1]_\mathbb{Q})\to\mathcal{L}_\Delta([0,1]_\mathbb{Q})$ by recursion as follows:
\begin{enumerate}
\item $t'(p)=p$ for $p\in Var$; $t'(\bot)=\bot$; $t'(\bar{c})=\bar{c}$;
\item $t'(\phi\circ\psi)=t'(\phi)\circ t'(\psi)$ for $\circ\in\{\land,\rightarrow\}$;
\item $t'(\Delta\phi)=\Delta t'(\phi)$;
\item $t'(\Box\phi)=\Box t'(\phi)$;
\item $t'([\phi]^1\bot)=\Delta t'(\phi)\rightarrow\bot$;
\item $t'([\phi]^1p)=\Delta t'(\phi)\rightarrow p$;
\item $t'([\phi]^1\bar{c})=\Delta t'(\phi)\rightarrow\bar{c}$;
\item $t'([\phi]^1(\psi\land\chi))=t'([\phi]^1\psi)\land t'([\phi]^1\chi)$;
\item $t'([\phi]^1(\psi\rightarrow\chi))=t'([\phi]^1\psi)\rightarrow t'([\phi]^1\chi)$;
\item $t'([\phi]^1\Box\psi)=\Delta t'(\phi)\rightarrow\Box t'([\phi]^1\psi)$;
\item $t'([\phi]^1\Delta\psi)=\Delta t'(\phi)\rightarrow\Delta t'([\phi]^1\psi)$;
\item $t'([\phi]^1[\psi]^1\chi)=t'([\phi\land[\phi]^1\psi]^1\chi)$.
\end{enumerate}
\end{definition}
As before, this translation provides equivalent formulas which contain no public announcement and this equivalence is provable, now in $\mathcal{GPA}^1_\Delta([0,1]_\mathbb{Q})$.
\begin{lemma}\label{lem:transval1}
For all $\phi\in\mathcal{L}_{PA}$: $\vdash_{\mathcal{GPA}^1_\Delta([0,1]_\mathbb{Q})}\phi\leftrightarrow t'(\phi)$.
\end{lemma}
The proof follows the standard procedure of defining a suitable complexity measure over which we then perform an induction as in Section \ref{sec:simplePAop}.
\begin{proof}
The key element is again a complexity measure $c:\mathcal{L}^\Delta_{PA^1}([0,1]_\mathbb{Q})\to\mathbb{N}$ on $\mathcal{L}_{PA}$ such that one can perform a suitable induction:
\begin{itemize}
\item $c(p)=c(\bot)=c(\bar{c})=1$;
\item $c(\phi\circ\psi)=1+\max\{c(\phi),c(\psi)\}$ for $\circ\in\{\land,\rightarrow\}$;
\item $c(\Box\phi)=c(\Delta\phi)=1+c(\phi)$;
\item $c([\phi]\psi)=(4+c(\phi))\cdot c(\psi)$.
\end{itemize}
This $c$ can be easily seen to have the following properties:
\begin{enumerate}
\item $c(\psi)\leq c(\phi)$ for $\psi$ a subformula of $\phi$;
\item $c(\Delta\phi\rightarrow\bot)<c([\phi]\bot)$;
\item $c(\Delta\phi\rightarrow p)<c([\phi]p)$;
\item $c(\Delta\phi\rightarrow\bar{c})<c([\phi]\bar{c})$;
\item $c([\phi]\psi\land[\phi]\chi)<c([\phi](\psi\land\chi))$;
\item $c([\phi]\psi\rightarrow[\phi]\chi)<c([\phi](\psi\rightarrow\chi))$;
\item $c(\Delta\phi\rightarrow\Box[\phi]\psi)<c([\phi]\Box\psi)$;
\item $c(\Delta\phi\rightarrow\Delta[\phi]\psi)<c([\phi]\Delta\psi)$;
\item $c([\phi\land[\phi]\psi]\chi)<c([\phi][\psi]\chi)$.
\end{enumerate}
Using the reduction axioms, it is now straightforward to prove the theorem by induction on $c(\phi)$.
\end{proof}
Again, this (combined with the internal definability of $\Delta$ via $[\cdot]^1$) yields the following expressivity result.
\begin{theorem}
$\mathcal{L}^\Delta_{PA^1}([0,1]_\mathbb{Q})\equiv\mathcal{L}^\Delta_\Box([0,1]_\mathbb{Q})\equiv\mathcal{L}_{PA^1}([0,1]_\mathbb{Q})$.
\end{theorem}
But, more importantly, we can even provide a completeness result for the previously introduced calculus.
\begin{theorem}
For any $\Gamma\cup\{\phi\}\subseteq\mathcal{L}^\Delta_{PA^1}([0,1]_\mathbb{Q})$, the following are equivalent:
\begin{enumerate}
\item $\Gamma\vdash_{\mathcal{GPA}^1_\Delta([0,1]_\mathbb{Q})}\phi$;
\item $\Gamma\models^\leq_{\mathsf{GK}}\phi$;
\item $\Gamma\models_{\mathsf{GK}}\phi$.
\end{enumerate}
\end{theorem}
The proof follows the same reasoning as Theorem \ref{thm:pacomp}.
\begin{remark}
As in the case without rational constants, the $\Delta$ operator is internally definable via $[\cdot]^1$. As discussed before, this makes it possible to obtain an axiomatization for $\mathcal{L}_{PA^1}([0,1]_\mathbb{Q})$ over all G\"odel-Kripke models from $\mathcal{GPA}^1_\Delta([0,1]_\mathbb{Q})$ by replacing the $\Delta$ with its translation. This yields the following proof calculus $\mathcal{GPA}^1([0,1]_\mathbb{Q})$ over the language $\mathcal{L}_{PA^1}([0,1]_\mathbb{Q})$ defined via:
\begin{description}
\item [($GK'$)] all the axiom schemes and rules of the calculus $\mathcal{GK}_\Delta([0,1]_\mathbb{Q})$ \emph{besides} ($\Delta 1$) - ($\Delta 5$), ($BK3$), ($\Delta\Box$) and ($N\Delta$);
\item [($\Delta 1'$)] $\neg[\phi]^1\bot\lor\neg\Delta\phi\neg[\phi]^1\bot$;
\item [($\Delta 2'$)] $\neg[\phi\lor\psi]^1\bot\rightarrow\neg[\phi]^1\bot\lor\neg[\psi]^1\bot$;
\item [($\Delta 3'$)] $\neg[\phi]^1\bot\rightarrow\phi$;
\item [($\Delta 4'$)] $\neg[\phi]^1\bot\rightarrow\neg[\neg[\phi]^1\bot]^1\bot$;
\item [($\Delta 5'$)] $\neg[\phi\rightarrow\psi]^1\bot\rightarrow(\neg[\phi]^1\bot\rightarrow\neg[\psi]^1\bot)$;
\item [($BK3'$)] $\neg[\bar{c}]^1\bot\leftrightarrow\bar{\delta(c)}$;
\item [($\Delta\Box'$)] $\neg[\Box\phi]^1\bot\rightarrow\Box\neg[\phi]^1\bot$;
\item [($N\Delta'$)] from $\phi$, infer $\neg[\phi]^1\bot$;
\item [($PA1'$)] $[\phi]^1\bot\leftrightarrow(\neg[\phi]^1\bot\rightarrow\bot)$;
\item [($PA2'$)] $[\phi]^1p\leftrightarrow(\neg[\phi]^1\bot\rightarrow p)$;
\item [($PA3'$)] $[\phi]^1\bar{c}\leftrightarrow(\neg[\phi]^1\bot\rightarrow \bar{c})$;
\item [($PA4'$)] $[\phi]^1(\psi\land\chi)\leftrightarrow([\phi]^1\psi\land[\phi]^1\chi)$;
\item [($PA5'$)] $[\phi]^1(\psi\rightarrow\chi)\leftrightarrow([\phi]^1\psi\rightarrow[\phi]^1\chi)$;
\item [($PA6'$)] $[\phi]^1\Box\psi\leftrightarrow(\neg[\phi]^1\bot\rightarrow\Box[\phi]^1\psi)$;
\item [($PA7'$)] $[\phi]^1\neg[\psi]^1\bot\leftrightarrow(\neg[\phi]^1\bot\rightarrow\neg[[\phi]^1]^1\bot)$;
\item [($PA8'$)] $[\phi]^1[\psi]^1\chi\leftrightarrow[\phi\land[\phi]^1\psi]^1\chi$.
\end{description}
Completeness for this system can be naturally proved by employing the translation which was used to defined it together with the above completeness result for $\mathcal{GPA}^1_\Delta([0,1]_\mathbb{Q})$.
\end{remark}
\section{The ugly: restrictive announcements}
Another public-announcement-type operator which we want to discuss is the restrictive announcement operator. To distinguish this new operator from the others, we consider yet again another language
\[
\mathcal{L}_{PA_r}:\phi::=\bot\mid p\mid (\phi\land\phi)\mid (\phi\rightarrow\phi)\mid\Box\phi\mid [\phi]^r\phi.
\]
Semantically, the language is again interpreted over G\"odel-Kripke models via the additional clause 
\[\vert [\phi]^r\psi\vert^w_\mathfrak{M}:=\vert\phi\vert^w_\mathfrak{M}\Rightarrow\vert\psi\vert^w_{\mathfrak{M}\vert^r_w\phi}
\]
for a G\"odel-Kripke model $\mathfrak{M}=\langle W,R,V\rangle$ with $w\in W$. But now, the model $\mathfrak{M}\vert^r_w\phi:=\langle W_{\mathfrak{M},w}^{\phi,r},R_{\mathfrak{M},w}^{\phi,r},V_{\mathfrak{M},w}^{\phi,r}\rangle$ is defined by restricting the set of worlds $W$ via
\[
W_{\mathfrak{M},w}^{\phi,r}:=\{v\in W\mid\vert\phi\vert^v_\mathfrak{M}\geq\vert\phi\vert^w_\mathfrak{M}\}
\] 
and $V_{\mathfrak{M},w}^{\phi,r}$ as well as $R_{\mathfrak{M},w}^{\phi,r}$ are just restrictions to this new set of worlds. Again, this extended evaluation function allows us to lift the consequence relation $\models_\mathsf{GK}$ to inputs $\Gamma\cup\{\phi\}\subseteq\mathcal{L}_{PA_r}$.\\

This operator may seem like another straightforward generalization of the classical semantics of public announcements, more so than the previous semantics given by modifying the accessibility function or announcing full truth as the usual way of defining public announcement in a classical semantic setting is by restricting the set of worlds relative to those which confine to the truth value of the formula in the world where the announcement is made. But, as the following lemma shows, this alternative generalization diverges drastically it is behavior compared to the other two.
\begin{lemma}
There is a $\mathsf{GK}$-model $\mathfrak{M}$ and a world $w\in\mathcal{D}(\mathfrak{M})$ such that
\[
\vert [\phi]^r\Box\psi\vert^w_\mathfrak{M}\neq\vert\phi\rightarrow\Box[\phi]^r\psi\vert^w_\mathfrak{M}.
\]
\end{lemma}
\begin{proof}
Consider the following model $\mathfrak{M}$
\begin{center}
\begin{tikzpicture}[shorten >=1pt,node distance=3cm,on grid,auto] 
   	\node[state] (q_0)  {$a$}; 
   	\node[state] (q_1) [right=of q_0] {$b$};
   	\path[<->] (q_0) edge  node {1} (q_1);
   	\path[->] 
    	(q_0) edge [loop above] node {1} ()
    	(q_1) edge [loop above] node {1} ();
\end{tikzpicture}
\end{center}
where $V(a,p)=V(a,q)=1$ and $V(b,p)=1/2$ but $V(b,q)=0$. Then
\begin{align*}
\vert[p]^r\Box q\vert^a_\mathfrak{M}&=V(a,p)\Rightarrow\vert\Box q\vert^a_{\mathfrak{M}\vert^r_ap}\\
&=R(a,a)\Rightarrow V(a,q)\\
&=1
\end{align*}
since $\mathfrak{M}\vert^r_ap$ only contains $a$ as a world. However, we have
\begin{align*}
\vert p\rightarrow\Box[p]^rq\vert^a_\mathfrak{M}&=V(a,p)\Rightarrow\vert\Box [p]^rq\vert^a_{\mathfrak{M}}\\
&=(R(a,a)\Rightarrow(V(a,p)\Rightarrow V(a,q)))\odot(R(a,b)\Rightarrow (V(b,p)\Rightarrow V(b,q)))\\
&=V(b,p)\Rightarrow V(b,q)\\
&=1/2\Rightarrow 0=0.
\end{align*}
\end{proof}
A similar problem occurs with the usual reduction axiom for iterated public announcement operators.
\begin{lemma}
There is a $\mathsf{GK}$-model $\mathfrak{M}$ and a world $w\in\mathcal{D}(\mathfrak{M})$ such that
\[
\vert[\phi]^r[\psi]^r\chi\vert^w_\mathfrak{M}\neq\vert[\phi\land[\phi]^r\psi]^r\chi\vert^w_\mathfrak{M}.
\]
\end{lemma}
\begin{proof}
Consider the following model $\mathfrak{M}$ with a similar frame as in the previous lemma:
\begin{center}
\begin{tikzpicture}[shorten >=1pt,node distance=3cm,on grid,auto] 
   	\node[state] (q_0)  {$a$}; 
   	\node[state] (q_1) [right=of q_0] {$b$};
   	\path[<->] (q_0) edge  node {1} (q_1);
   	\path[->] 
    	(q_0) edge [loop above] node {1} ()
    	(q_1) edge [loop above] node {1} ();
\end{tikzpicture}
\end{center}
However, the evaluations change to the following: we define
\begin{enumerate}
\item $V(a,p):=1$, $V(a,q):=1/2$, $V(a,s):=1$,
\item $V(b,p):=1/2$, $V(b,q):=1$, $V(b,s):=0$.
\end{enumerate}
Then, naturally, $\mathfrak{M}\vert^r_a p$ is just $\mathfrak{M}$ restricted to $\{a\}$, which we denote here by $\mathfrak{M}\upharpoonright\{a\}$, and therefore also $(\mathfrak{M}\vert^r_a p)\vert^r_a q=\mathfrak{M}\upharpoonright\{a\}$. This yields
\begin{align*}
\vert[p]^r[q]^r\Box s\vert^a_\mathfrak{M}&=V(a,p)\Rightarrow\vert[q]^r\Box s\vert^a_{\mathfrak{M}\vert^r_a p}\\
&=V(a,p)\Rightarrow\left(V(a,q)\Rightarrow\vert\Box s\vert^a_{(\mathfrak{M}\vert^r_a p)\vert^r_a q}\right)\\
&=1\Rightarrow\left(1/2\Rightarrow\vert\Box s\vert^a_{\mathfrak{M}\upharpoonright\{a\}}\right)\\
&=1/2\Rightarrow V(a,s)=1.\\
\end{align*}
On the other hand, we have
\begin{align*}
\vert p\land[p]^rq\vert^a_\mathfrak{M}&=V(a,p)\odot\left(V(a,p)\Rightarrow\vert q\vert^a_{\mathfrak{M}\vert^r_a p}\right)\\
&=1\odot\left(1\Rightarrow V(a,q)\right)=1/2
\end{align*}
and similarly
\begin{align*}
\vert p\land[p]^rq\vert^b_\mathfrak{M}&=V(b,p)\odot\left(V(b,p)\Rightarrow\vert q\vert^b_{\mathfrak{M}\vert^r_b p}\right)\\
&=1/2\odot\left(1/2\Rightarrow V(a,q)\right)=1/2.
\end{align*}
Therefore $\vert p\land[p]^rq\vert^b_\mathfrak{M}\geq\vert p\land[p]^rq\vert^a_\mathfrak{M}$ and thus $\mathfrak{M}\vert^r_a(p\land [p]^rq)=\mathfrak{M}$. But this yields
\begin{align*}
\vert[p\land [p]^rq]^r\Box s\vert^a_\mathfrak{M}&=\vert p\land[p]^rq\vert^a_\mathfrak{M}\Rightarrow\vert\Box s\vert^a_{\mathfrak{M}\vert^r_a(p\land [p]^rq)}\\
&=1/2\Rightarrow\vert\Box s\vert^a_\mathfrak{M}\\
&=1/2\Rightarrow 0=0.
\end{align*}
\end{proof}
We leave it at these observations for this note but want to remark that it does not seem too far fetched to pursue a reduction-style completeness proof to a (suitable extension of) $\mathcal{GK}$ as the $[\cdot]^r$-operator does commute with the basic propositional connectives.
\begin{lemma}
The following formulas are valid in all $\mathsf{GK}$-models:
\begin{enumerate}
\item $[\phi]^r\bot\leftrightarrow(\phi\rightarrow\bot)$;
\item $[\phi]^rp\leftrightarrow(\phi\rightarrow p)$;
\item $[\phi]^r(\psi\land\chi)\leftrightarrow([\phi]^r\psi\land[\phi]^r\chi)$;
\item $[\phi]^r(\psi\rightarrow\chi)\leftrightarrow([\phi]^r\psi\rightarrow[\phi]^r\chi)$.
\end{enumerate}
\end{lemma}
The proofs follows the same vein as the ones of Lemma \ref{lem:fpavalidity}.\\

We thus leave it at the following question:
\begin{question}
Is $\mathcal{L}_{PA_r}\equiv\mathcal{L}_\Box$ or do we have $\mathcal{L}_\Box\prec\mathcal{L}_{PA_r}$?
\end{question}
\bibliographystyle{plain}
\bibliography{ref}
\end{document}